\documentclass[a4paper,10pt]{amsart}
\usepackage{graphicx} 
\usepackage[a4paper, margin=2 cm]{geometry}
\usepackage{hyperref}
\usepackage{amsthm}
\usepackage{thm-restate}
\usepackage{geometry}
\usepackage{amssymb}
\usepackage{amsmath}
\usepackage{cleveref}
\usepackage{comment}
\usepackage{setspace}
\usepackage{xcolor}
\usepackage[shortlabels]{enumitem}
\usepackage{listings}
\usepackage{blkarray}
\usepackage{lmodern}
\usepackage[english]{babel}
\usepackage{float}
\usepackage{tikz}
\usepackage{mathtools}
\usetikzlibrary{arrows.meta,arrows}
\usetikzlibrary{arrows,decorations.markings}
\usetikzlibrary{decorations.pathmorphing}

\newtheorem{theorem}{Theorem}[section]
\newtheorem{lemma}[theorem]{Lemma}

\newtheorem{corollary}[theorem]{Corollary}

\newtheorem*{claim*}{Claim}

\theoremstyle{definition}
\newtheorem{definition}[theorem]{Definition}

\title{Improved bounds for zero-sum cycles in $\mathbb{Z}_p^d$}
\author
{
Micha Christoph
}
\author 
{
Charlotte Knierim
}
\author
{
Anders Martinsson
}
\author
{
Raphael Steiner 
}
\address[Christoph,Knierim,Martinsson,Steiner]{Department of Computer Science, Institute of Theoretical Computer Science, ETH Z\"{u}rich, Switzerland}
\email{\tt $\{$micha.christoph,charlotte.knierim,anders.martinsson,raphaelmario.steiner$\}$@inf.ethz.ch}

\thanks{The first and fourth authors were supported by the SNSF Ambizione Grant No. 216071.}

\date{\today}

\begin{document}
\begin{abstract}
For a finite  Abelian group $(\Gamma,+)$, let $n(\Gamma)$ denote the smallest positive integer $n$ such that for each labelling of the arcs of the complete digraph of order $n$ using elements from $\Gamma$, there exists a directed cycle such that the total sum of the arc-labels along the cycle equals $0$. Alon and Krivelevich initiated the study of the parameter $n(\cdot)$ on cyclic groups and proved that $n(\mathbb{Z}_q)=O(q\log q)$. Several improvements and generalizations of this bound have since been obtained, and an optimal bound in terms of the group order of the form $n(\Gamma)\le |\Gamma|+1$ was recently announced by Campbell, Gollin, Hendrey and the last author. While this bound is tight when the group $\Gamma$ is cyclic, in cases when $\Gamma$ is far from being cyclic, significant improvements on the bound can be made. In this direction, studying the prototypical case when $\Gamma=\mathbb{Z}_p^d$ is a power of a cyclic group of prime order, Letzter and Morrison showed that $n(\mathbb{Z}_p^d) \le O(pd(\log d)^2)$ and that $n(\mathbb{Z}_2^d)\le O(d \log d)$. They then posed the problem of proving an (asymptotically optimal) upper bound of $n(\mathbb{Z}_p^d)\le O(pd)$ for all primes $p$ and $d \in \mathbb{N}$. In this paper, we solve this problem for $p=2$ and improve their bound for all primes $p \ge 3$ by proving $n(\mathbb{Z}_2^d)\le 5d$ and $n(\mathbb{Z}_p^d)\le O(pd\log d)$. While the first bound determines $n(\mathbb{Z}_2^d)$ up to a multiplicative error of $5$, the second bound is tight up to a $\log d$ factor. Moreover, our result shows that a tight bound of $n(\mathbb{Z}_p^d)=\Theta(pd)$ for arbitrary $p$ and $d$ would follow from a (strong form) of the well-known conjecture of Jaeger, Linial, Payan and Tarsi on additive bases in $\mathbb{Z}_p^d$. 

Along the way to proving these results, we establish a generalization of a hypergraph matching result by Haxell in a matroidal setting. Concretely, we obtain sufficient conditions for the existence of matchings in a hypergraph whose hyperedges are labelled by the elements of a matroid, with the property that the edges in the matching induce a basis of the matroid. We believe that these statements are of independent interest.
\end{abstract}
\maketitle
\section{Introduction}
\emph{Zero-sum Ramsey theory} forms a large branch of modern Ramsey theory that lies at the intersection of combinatorics, number theory, algebra, graph theory, discrete analysis and other branches of mathematics. Problems in this area are typically concerned with combinatorial objects that are equipped with labels from a given group, and study conditions on these objects which guarantee the containment of a substructure that has a zero total sum in the group associated with the labels (this is in contrast but analogous to the usual Ramsey-type problems, in which one seeks to find monochromatic substructures). The history of such problems dates back at least to work of Erd\H{o}s, Ginzburg and Ziv~\cite{EGZ} from 1961, who proved that every sequence of $2m-1$ elements in the cyclic group $(\mathbb{Z}_m,+)$ contains a subsequence of length $m$ of total sum equal to zero. Since this fundamental result, many more zero-sum Ramsey problems have been studied, covering a variety of combinatorial structures such as graphs, digraphs and hypergraphs. For more background on such problems, we refer to Caro's survey~\cite{C} on zero-sum Ramsey theory. 

In the paper at hand, we shall be concerned with a zero-sum Ramsey problem for directed cycles in complete digraphs\footnote{A \emph{complete digraph} is a digraph which includes both arcs $(x,y), (y,x)$, for every two distinct vertices $x$ and $y$ in its vertex-set.}, that was first studied by Alon and Krivelevich~\cite{AK}. Let~$D$ be a digraph and let~${(\Gamma,+)}$ be an Abelian group. 
A map~${w \colon A(D) \to \Gamma}$ (where $A(D)$ denotes the arc set of $D$) is called a \emph{$\Gamma$-labelling of~$D$}, 
and the pair~${(D,w)}$ is called a \emph{$\Gamma$-labelled digraph}. 
We say that a directed cycle~$C$ in $D$ with vertex set ${\{ v_i : i \in \mathbb{Z}_\ell \}}$ and arc set~${\{ (v_{i}, v_{i+1}) : i \in \mathbb{Z}_\ell \}}$ is \emph{$w$-zero} if~${\sum_{i=0}^{\ell-1}{w(v_{i}, v_{i+1})}}$ equals the neutral element $0$ of $\Gamma$. When $w$ is clear from context, we simply call $C$ a \emph{zero-sum cycle}. The zero-sum Ramsey parameter we study in this paper is denoted by $n(\Gamma)$ and defined as the smallest positive integer $n$ such that every $\Gamma$-labelled complete digraph on $n$ vertices contains a zero-sum cycle. Alon and Krivelevich~\cite{AK} initiated the study of $n(\Gamma)$ for cyclic groups and proved the bound $n(\mathbb{Z}_q)\le O(q \log q)$ via a beautiful probabilistic argument. They also showed an improved bound of $n(\mathbb{Z}_p)\le 2p-1$ when $p$ is a prime. They then used these bounds to deduce the existence of cycles of length divisible by $q$ in all graphs of sufficiently large Hadwiger number in terms of $q$. More precisely, they showed that every graph containing a complete graph on $2n(\mathbb{Z}_q)$ vertices as a minor (i.e., of Hadwiger number at least $2n(\mathbb{Z}_q)$) necessarily contains a cycle of length divisible by $q$. More generally, they proved the existence of so-called \emph{divisible subdivisions} in graphs of large Hadwiger number, see also the recent paper~\cite{DDS} which provides asymptotically optimal bounds for this problem. 

The result of Alon and Krivelevich on cyclic groups was improved to a linear bound $n(\Gamma)\le 8|\Gamma|$ for all Abelian groups $(\Gamma,+)$ by M\'{e}szar\'{o}s and the last author~\cite{MS}, who also showed that $n(\mathbb{Z}_p)\le \frac{3p}{2}$ for every prime number $p$. A further improvement of the general bound was achieved independently by Berendsohn, Boyadzhiyska and Kozma~\cite{BBK} and by Akrami et al.~\cite{AACGMM}, who showed that $n(\Gamma)\le 2|\Gamma|-1$ for every non-trivial $\Gamma$. Campbell, Gollin, Hendrey and the last author~\cite{CGHS} recently announced an optimal bound in terms of the group order of the form $n(\Gamma)\le |\Gamma|+1$. While one can see that this latter bound is tight for cyclic groups (see e.g.~\cite{MS}), it does not seem to be tight for many groups that are ``far'' from being cyclic. Letzter and Morrison~\cite{ML} recently provided evidence in this direction by studying the prototypcial case of bounding $n(\Gamma)$ for large powers of cyclic groups of prime order. As their main results, they proved that $n(\mathbb{Z}_p^d)\le O(pd(\log d)^2)$ for every prime number $p$ and every $d \in \mathbb{N}$, with an improved bound of $n(\mathbb{Z}_2^d)\le O(d \log d)$ for the case $p=2$. They then posed the problem of improving these bounds, and in particular raised the question whether $n(\mathbb{Z}_p^d)\le O(pd)$ for every prime number $p$ and $d \in \mathbb{N}$. 

\medskip

\paragraph*{\textbf{Our results.}}
In this paper, we make substantial progress on the aforementioned problem of Letzter and Morrison, by affirmatively solving their problem in the case $p=2$ and asymptotically improving on their bound for the case of general $p$. To state our result, we need the following notion from~\cite{ALM}: For $p$ a prime number and $d \in \mathbb{N}$, $f(p,d)$ denotes the smallest positive integer satisfying the following property: For every collection $B_1,\ldots,B_{f(p,d)}$ of linear bases of the vector-field $\mathbb{Z}_p^d$ over the prime-field $\mathbb{Z}_p$, the multi-set union $B_1\cup \dots \cup B_{f(p,d)}$ forms an additive basis\footnote{Recall that given a vector space $S$, an \emph{additive basis} of $S$ is a (multi-)set $B$ consisting of elements from $S$ such that every member of $S$ can be written as the sum of elements in a (possibly empty) (multi-)subset of $B$.} of $\mathbb{Z}_p^d$. Alon, Linial and Meshulam~\cite{ALM} proved the (still currently best) bound of $f(p,d)\le (p-1)\log d + (p-2)$. A central problem in additive combinatorics connected to this parameter is the so-called \emph{additive basis conjecture}, which states that $f(p,d)$ can be upper-bounded by a function of $p$ only~\cite{ALM}. 
Moreover, that $f(p,d)\le cp$ for some absolute constant $c>0$ and all primes $p$ and $d \in \mathbb{N}$, or even that $f(p,d)=p$, seems plausible~\cite{OP}. Using this definition, we can state the main result of our paper as follows.

\begin{theorem}\label{thm:binary}
For every positive integer $d$, we have $n(\mathbb{Z}_p^d)\le 5d\cdot f(p,d)$.
\end{theorem}
Using the above-mentioned bound on $f(p,d)$ by Alon, Linial and Meshulam~\cite{ALM} as well as the fact that $f(2,d)=1$ for every\footnote{This follows since the notions of linear and additive bases coincide in the case of $\mathbb{Z}_2^d$.} $d$ we have the following immediate consequences.
\begin{corollary}
For every prime $p\ge 3$ and $d \in \mathbb{N}$, we have $n(\mathbb{Z}_p^d)= O(pd\log d)$.
\end{corollary}
\begin{corollary}
    For every $d \in \mathbb{N}$, we have $n(\mathbb{Z}_2^d)\le 5d$. 
\end{corollary}

We note that the above corollary implies that $n(\mathbb{Z}_p^d)=O(pd)$ would follow from $f(p,d)=O(p)$ (i.e., a strong quantitative form of the additive basis conjecture), and thus answer the question of Letzer and Morrison affirmatively for all values of $p$ and $d$. 
As already noted in~\cite{ML}, this bound would be asymptotically optimal, as it is easy to observe that $n(\mathbb{Z}_p^d)> (p-1)d$. In particular, our upper bound $n(\mathbb{Z}_2^d)\le 5d$ is tight up to a multiplicative constant of at most $5$.  

As a key technical lemma to show these results, we prove a generalization of Haxell's condition for the existence of perfect matchings in bipartite hypergraphs \cite{Haxell95}. Haxell's condition can equivalently be stated as follows. Given any hypergraph $H$ in which the size of every hyperedge is at most $r\geq 1$, and a coloring of the hyperedges of $H$ with $d$ colors, one of the following holds:
\begin{enumerate}
    \item $H$ contains a rainbow matching, i.e.\ a matching with one hyperedge of each color, or
    \item there exists a $k\geq 1$ and a vertex set $U$ of size at most $(2r-1)(k-1)$ such that at most $d-k$ colors are present in $H-U.$
\end{enumerate}
Note that in the latter case, we may assume that $H-U$ contains a rainbow matching, as otherwise we can simply repeatedly apply the statement to $H-U$ until we reach a subhypergraph where this is the case. Here we consider a generalization of this setting, where the hyperedges of $H$ are labelled with elements from a matroid $\mathcal{M}$ and where the notion of rainbow matching is replaced by a matching whose labels form a basis of $\mathcal{M}$ (note that this contains rainbow matchings as a special case, namely when $\mathcal{M}$ is the free matroid on the set of colors). The direct analogue of Haxell's condition is given in Lemma \ref{lemma:matroidtransversal}. Given the well-known generalization of bipartite matching to matroid intersection, as first shown by Edmonds \cite{Edmonds70}, see also Chapter 41 of \cite{Schrijver}, this seems like a natural generalization of bipartite hypergraph matching. Despite intensive search, we could not find these type of statements in the literature, and to the best of our knowledge they are novel. We refer to the works of Aharoni and Berger~\cite{AB} and Davies, Rothvoss and Zhang~\cite{DRZ} for some results that are similar in spirit, but do not seem directly related.  We think that our generalization may be of independent interest. In particular, we expect that many more applications of our statements on independent hypergraph matchings to combinatorial and zero-sum problems will be found in the future.

\section{Notation}
A (multi)-hypergraph is a tuple $H=(V,E)$, where $V$ is a finite set and $E$ is a multi-set whose members are non-empty subsets of $V$. For a (multi-)hypergraph $H$, we use the notation $V(H)$ and $E(H)$ to refer to its vertex-set and hyperedge-(multi-)set. As usual, for a subset $U\subseteq V(H)$ of vertices we denote by $H-U$ the (multi-)hypergraph with vertex-set $V(H)\setminus U$ that inherits from $H$ all hyperedges that are disjoint from $U$.
With a slight abuse of notation, for a (multi-)subset $X\subseteq E(H)$ we write $V(X):=\bigcup_{e\in X}{e}$ for the set of vertices covered by the hyperedges in $X$. Finally, we say that a (multi-)hypergraph $H$ is connected if for every partition of its vertex-set into non-empty parts $A, B$ there exists a hyperedge containing a vertex from both $A$ and $B$. 

Throughout the rest of the paper, we use standard matroid notation, as can be found e.g.\ in the classic book~\cite{O} by Oxley. Specifically, given a matroid $\mathcal{M}$, we denote by $E(\mathcal{M})$ its set of elements, and for a subset $X\subseteq E(\mathcal{M})$ by $\mathrm{span}(X)$ the span of $X$ in the matroid $\mathcal{M}$. We will also need the standard definition of \emph{direct sums} of matroids. Given a sequence of matroids $\mathcal{M}_1,\ldots,\mathcal{M}_k$, their \emph{direct sum} $\mathcal{M}_1+\cdots+\mathcal{M}_k$ is defined as the matroid with element-set equal to the disjoint union $E(\mathcal{M}_1)\cup\cdots\cup E(\mathcal{M}_k)$ in which a subset $I$ of elements is independent if and only if $I \cap E(\mathcal{M}_i)$ is independent in $\mathcal{M}_i$, for every $i \in \{1,\ldots,k\}$. It is well-known and easy to verify that this operation indeed always results in a matroid again.

Our arguments and statements make use of the notion of an independent matching in a hypergraph equipped with matroid-labels, formally defined as follows.

\begin{definition}[independent matching]
    Let $H$ be a (multi-)hypergraph, $\mathcal{M}$ a matroid and $\gamma\colon E(H)\rightarrow E(\mathcal{M})$ a hyperedge-labelling\footnote{Note here that we allow $\gamma$ to assign distinct matroid elements to distinct (multi-)hyperedges of $H$ independently, even if they form the same subset of vertices in $H$.} of $H$. An \emph{independent matching} $M$ in $(H,\gamma)$ is a collection of pairwise vertex disjoint hyperedges in $H$ with distinct hyperedge-labels such that $\gamma(M)$ forms an independent set in $\mathcal{M}$. 

    We say that an independent matching $M$ in $(H,\gamma)$ is \emph{inclusion-wise maximal} if there exists no independent matching $M'$ such that $M\subsetneq M'$. We say that $M$ is \emph{maximal}, if there exists no independent matching $M'$ in $(H,\gamma)$ such that $\mathrm{span}(\gamma(M')) = \mathrm{span}(\gamma(M))$ and $M'$ is not inclusion-wise maximal. 
\end{definition}
Note that, importantly, the matching consisting of an empty set of hyperedges will also be independent with our above definition.

\section{Proof of the main result}

Before we dive into the details of the proof, we establish facts about independent matchings in hypergraphs. The proof method of the following first result is in spirit similar to the classic augmenting path arguments for matchings in simple graphs.

\begin{lemma}\label{lemma:key}
    Let $H$ be a (multi-)hypergraph, $\mathcal{M}$ a matroid and $\gamma\colon E(H)\rightarrow E(\mathcal{M})$.
    Let $M$ be a maximal independent matching of size $\ell$ in $H$. 
    Suppose there exists $e\in E(H)$ such that $\gamma(e)\notin \mathrm{span}(\gamma(M))$. Then, there exists $k>0$ and a collection of $k+1$ hyperedges $X\subseteq E(H)$ such that $H-V(X)$ contains a maximal independent matching of size $\ell-k$ and the subhypergraph of $H$ formed by the hyperedges in $X$ and their vertices
    is connected.
\end{lemma}
\begin{proof} In the following proof, we write $e\sqcap M$ to denote the set of hyperedges in a matching $M$ that intersect a hyperedge $e$. 

    Let $M$ be a maximal independent matching in $H$ of size $\ell$ and let $S=\mathrm{span}(\gamma(M))$. Take $e\in E(H)$ with $\gamma(e)\notin S$ given by assumption. Without loss of generality, let us assume that $M$ is an independent matching that minimizes $|e\sqcap M|$ subject to $\mathrm{span}(\gamma(M))=S$.  
    
    Observe that, as $M$ is inclusion-wise maximal, $e\sqcap M$ is non-empty. 
    Let $X=\{e\}\cup (e\sqcap M)$ and let $k>0$ denote the number of hyperedges in $e\sqcap M$. Hence, $X$ contains $k+1$ hyperedges and is connected, since $e$ intersects all the hyperedges in $X$. Additionally, $M\setminus X$ is an independent matching in $H-V(X)$ of size $\ell-k$. To conclude the proof, we show that $M\setminus X$ is maximal in $H-V(X)$.
    
    Let $S'=\mathrm{span}(\gamma(M\setminus X))$. Suppose towards a contradiction that there exists an independent matching $M'$ in $H-V(X)$ that spans $S'$ but is not inclusion-wise maximal. Let $e^*\in E(H-V(X))$ be a hyperedge that extends $M'$ to a larger independent matching in $H-V(X)$. Let $M''=M'\cup (M\cap X)$. Since $\gamma(M')$ spans $S'$, we have that $M''$ is an independent matching with $\mathrm{span}(\gamma(M''))=S$. Therefore, $M''$ is inclusion-wise maximal in $H$. Since $M''$ is inclusion-wise maximal, we get $\gamma(e^*)\in S$. Then, $\gamma(e^*)\in\mathrm{span}(\gamma(M''))\setminus\mathrm{span}(\gamma(M'))$. Thus, there exists $e'\in M''\setminus M'=M\cap X$ such that $e^*$ is not in the span of $\gamma(M''\setminus \{e'\})$. It follows that $M^*=M''\cup\{e^*\}\setminus\{e'\}$ is an independent matching with span $S$. But then, $|e\sqcap M^*|=|(e\sqcap M)\setminus \{e'\}|=|e\sqcap M|-1$ contradicting the minimality of $|e\sqcap M|$.
    \end{proof}

With this tool at hand, we are ready to prove our main technical lemma.

\begin{lemma}\label{lemma:matroidtransversal} Let $H$ be a (multi-)hypergraph in which the size of every hyperedge is at most $r$, let $\mathcal{M}$ be a matroid of rank $d\geq 1$ and $\gamma\colon E(H)\rightarrow E(\mathcal{M})$. Then, one of the following statements hold:
\begin{enumerate}
    \item $H$ contains an independent matching of size $d$, or
    \item there exists a $k\geq 1$ and a vertex set $U$ of size at most $(2r-1)(k-1)$ such that $H-U$ contains an independent matching $M$ of size $d-k$ with $\mathrm{span}(\gamma(M))=\mathrm{span}(\gamma(E(H-U))).$
\end{enumerate}
\end{lemma}
\begin{proof} 
We prove the lemma by considering the following slightly stronger statement: For any $(H,\gamma)$ as above that contains a maximal independent matching of size $\ell\geq 0$, there exists $0\leq a\leq \ell$ and $U\subseteq V(H)$ of size at most $(2r-1)a$ such that $H-U$ contains an independent matching of size $\ell-a$ whose labels form a basis of $\mathrm{span}(\gamma(E(H-U))).$

To see that this implies the lemma, let $M$ be any maximal independent matching of $H$. If $|M|=d$, then conclusion $(1)$ holds. Otherwise $\ell=|M|\leq d-1$, and there exists $0\leq a \leq \ell$ and a vertex set $U$ of size at most $(2r-1)a$ such that $H-U$ contains an independent matching of size $\ell-a$ whose labels form a basis of $\mathrm{span}(\gamma(E(H-U))).$ One can immediately check that conclusion $(2)$ of the lemma holds for $k=d-\ell+a$. Note that, as $\ell\leq d-1$, it follows that $k\geq a+1$, which implies $k\geq 1$ and $|U|\leq (2r-1)(k-1)$ as desired.

We now prove the above statement by induction on $\ell$. Let $M$ be a maximal independent matching of $H$ of size $\ell\geq 0$. If $\mathrm{span}(\gamma(M))=\mathrm{span}(\gamma(E(H)))$, then the induction hypothesis follows with $a=0$ and $U=\emptyset$. Note that this resolves the case when $\ell = 0$. Otherwise, $\ell\leq d-1$ and there exists $e\in E(H)$ such that $\gamma(e)\not\in \mathrm{span}(\gamma(M)).$ By Lemma \ref{lemma:key} there exists $0< k \leq \ell$ and $X\subseteq E(H)$ with $|X|=k+1$ such that the hyperedges in $X$ form a connected subhypergraph of $H$ and $H-V(X)$ contains a maximal independent matching of size $\ell-k.$ By the induction hypothesis, there exists an $0\leq a' \leq \ell-k$ and a vertex set $U'$ of size at most $(2r-1)a'$ such that $H-V(X)-U'$ contains an independent matching of size $\ell-k-a'$ whose labels form a basis of $\mathrm{span}(\gamma(E(H-V(X)-U')))$. We claim that the statement holds with $U=V(X)\cup U'$ and $a=k+a'$. Indeed, $a\geq 0$ and $H-U$ contains an independent matching of size $\ell-a$ whose labels form a basis of $\mathrm{span}(\gamma(E(H-U))),$ so it only remains to check that $|U|\leq (2r-1)a.$ Since $|X|$ is connected, we get 
$$|U|= |V(X)|+|U'|\leq (r-1)k + r + (2r-1)a' = (2r-1)(k+a') - r(k-1) \leq (2r-1)a. $$
The statement follows by induction.
\end{proof}
The following lemma tailors \Cref{lemma:matroidtransversal} to our use case.
\begin{lemma}\label{lemma:differentbases} Let $H$ be a (multi-)hypergraph in which the size of every hyperedge is at most $r$, let $\mathcal{M}$ be a matroid of rank $d\geq 1$ and $\gamma\colon E(H)\rightarrow E(\mathcal{M})$. For any $m\geq 1$ there exists $U\subseteq V(H)$ of size at most $(2r-1)(dm-1)$ and a sequence $M_1, \dots, M_m$ of $m$ pairwise vertex disjoint matchings in $H-U$ such that, for each $1\leq i \leq m$, $\gamma(M_i)$ is a basis of $\mathrm{span}(\gamma(E(H-U))).$ 
\end{lemma}
\begin{proof} 
Given $(H,\gamma)$ as above we form a new hypergraph $H'$ on the same vertex set with hyperedge-labels from the $m$-fold direct sum $\mathcal{M}'=\mathcal{M}+\dots + \mathcal{M}$. Note that $\mathcal{M}'$ has rank $dm$. For each hyperedge $e$ in $H$ we place $m$ copies of it in $H'$ and label each of these with a distinct copy of $\gamma(e)$ from the summands in $\mathcal{M}'$.

We claim there exists a vertex set $U$ of size at most $(2r-1)(dm-1)$ such that $H'-U$ contains an independent matching whose labels form a basis of $\mathrm{span}(\gamma(E(H'-U)))$. Indeed, applying Lemma \ref{lemma:matroidtransversal} to $H'$ we either get from case $(1)$ that the statement holds with $U=\emptyset$, or from case $(2)$ with the set $U$ as given in the lemma. By construction of $H'$ we have $\mathrm{span}(\gamma(E(H'-U)))=\mathrm{span}(\gamma(E(H-U)))+\dots+\mathrm{span}(\gamma(E(H-U)))$. Hence, letting $M_i$ be the set of hyperedges in $H$ whose $i$-th copy in $H'$ is in $M$ for $1\leq i \leq m$ we obtain $m$ matchings of $H-U$ with the desired properties.
\end{proof}

We are now ready to apply Lemma~\ref{lemma:differentbases} to deduce our main result, Theorem~\ref{thm:binary}.

\begin{proof}[Proof of Theorem~\ref{thm:binary}]
Let $D$ be a complete digraph of order $5d\cdot f(p,d)$ and let $w:A(D)\rightarrow \mathbb{Z}_p^d$ be an arc-labelling of $D$. We need to show that $D$ contains a directed cycle such that its arc-labels sum to zero. To achieve this goal, consider the $3$-uniform multi-hypergraph $H$ that has vertex-set $V(D)$ and one hyperedge for each ordered triple $(x, y, z)$ of distinct vertices, associated with the vertex set $\{x, y, z\}$. In particular, every $3$-subset of $V(D)$ is included in $E(H)$ in exactly $6$ different hyperedges.

Furthermore, define a labelling $\gamma:E(H)\rightarrow \mathbb{Z}_p^d$ as follows. For every hyperedge $(x, y, z)$ in $H$, we define $\gamma(x, y, z):=w(x,y)+w(y,z)-w(x,z)$. Note that $\mathbb{Z}_p^d$ can be viewed as a $d$-dimensional vector space over the prime field $\mathbb{F}_p$ and as such gives rise to a linear matroid $\mathcal{M}$ with element-set $E(\mathcal{M})=\mathbb{Z}_p^d$. 

Now, apply Lemma~\ref{lemma:differentbases} to the pair $(H,\mathcal{M})$ and the labelling $\gamma$, with parameter $m:=f(p,d)$. Accordingly, we obtain that there exists $U \subseteq V(D)$ with $|U|\le (2\cdot 3-1)\cdot (d\cdot f(p,d)-1)=5d\cdot f(p,d)-5$ and a sequence $M_1,\ldots,M_{f(p,d)}$ of pairwise vertex-disjoint matchings in $H-U$ such that, for each $1 \le i \le f(p,d)$, we have that $\gamma(M_i)$ is a basis of $\mathrm{span}(\gamma(E(H-U))$. Suppose for a moment that $\mathrm{span}(\gamma(E(H-U))=\{0\}$. Let $x,y,z$ be any $3$ distinct vertices in $V(D)\setminus U$. (Recall, by assumption of $D$ and $U$, $|V(D)\setminus U|\geq 5$). Then $\gamma(x, y, z)+\gamma(x, z, y) = w(y, z)+w(z, y)=0$. Hence, in this case we have found a zero-sum directed cycle of length two, concluding the proof. So, in the following, assume that $\mathrm{span}(\gamma(E(H-U))\neq \{0\}$.

Next, note that $\mathrm{span}(\gamma(E(H-U))$ is a vector space isomorphic to $\mathbb{Z}_p^{d'}$ for some $1\le d'\le d$, and since we have $f(p,d)\ge f(p,d')$, it follows by definition that $\gamma(M_1\cup \cdots \cup M_{f(p,d)})$ forms an additive basis of $\mathrm{span}(\gamma(E(H-U))$. Let us enumerate the hyperedges in $M:=M_1\cup \cdots \cup M_{f(p,d)}$ as $(x_i,y_i,z_i), i=1,\ldots,|M|$. Now consider the collection $\mathcal{Q}$ of directed cycles 
in $D$ that are obtained from the ``base cycle'' $Q_0=x_1,z_1,\ldots,x_{|M|},z_{|M|},x_1$ by replacing an arbitrary subset of the arcs $(x_i,z_i)$ with the length two-segments $(x_i,y_i,z_i)$. Denoting by $a:=\sum_{e \in A(Q_0)}{w(e)}$ the sum of labels in $Q_0$, we can then see that for every element $x\in a+\mathrm{span}(\gamma(E(H-U))$, there exists a cycle $Q \in \mathcal{Q}$ such that $\sum_{e \in A(Q)}{w(e)}=x$. We now claim that $a \in \mathrm{span}(\gamma(E(H-U))$. To see this, fix any vertex $u \notin V(Q_0)$, for instance $u=y_1$, and note that $$a = \sum_{(v_1,v_2) \in A(Q_0)}{(w(u,v_1)+w(v_1,v_2)-w(u,v_2))},$$ which, as a sum of $\gamma$-labels of hyperedges in $D-U$, belongs to $\mathrm{span}(\gamma(E(H-U))$. Finally, this shows that, for any vector in $\mathrm{span}(\gamma(E(H-U)))$, there exists a cycle $Q\in\mathcal{Q}$ whose edge labels sum up to this vector. In particular, $\mathcal{Q}$ contains a zero-sum cycle as desired.

\end{proof}

\end{document}